\documentclass[12pt,a4paper]{amsart}
\usepackage[english]{babel}
\usepackage[utf8]{inputenc}
\usepackage{comment}
\usepackage{amssymb}
\usepackage{amsmath} 
\usepackage{amsthm} 
\usepackage{enumitem}
\usepackage{bbm}
\usepackage{bbold}
\usepackage{todonotes}
\usepackage{verbatim}

\newtheorem{thm}{Theorem}[section]
\newtheorem{re}[thm]{Remark}

\newtheorem{prop}[thm]{Proposition}
\newtheorem{lemma}[thm]{Lemma}
\newtheorem{cor}[thm]{Corollary}
\numberwithin{equation}{section}

\DeclareMathOperator*{\res}{ res}

\title{A Voronoi summation formula for the shifted triple divisor function}
\author{Alessandro Fazzari}
\address{American Institute of Mathematics, 600 East Brokaw Road, San Jose, CA 95112, US}
\email{fazzari@aimath.org}

\begin{document}

\maketitle

\begin{abstract}
In this paper, we prove a Voronoi summation formula for the shifted 3-fold divisor function twisted by additive characters. As the main tool, we provide the functional equation for the shifted $GL(3)$ Estermann function.
\end{abstract}

\section{Introduction and statement of the main result}

Voronoi summation formulas play a pivotal role in number theory and they have been used intensively for many applications in the past years (see \cite{MS-survey} for a historical detour about summation formulas). They can be seen as Poisson formulas weighted by Fourier coefficients of automorphic forms and used to study analytic properties of the attached $L$-functions. 
A classical example is given by the Dirichlet's divisor problem, involving the mean value of the divisor function $\tau(n)$. Dirichlet showed that
$$\sum_{n\leq x}\tau(n) -x\log x - (2\gamma-1)x = O(\sqrt{x}) , $$
where $\gamma$ is the Euler's constant. Voronoi's improvement on Dirichlet's theorem is due to a summation formula for the divisor function, which can be stated as follows
\begin{equation}\begin{split}\label{VoronoiSummationFormula}  
\sum_{n\leq x} \tau(n)\phi(n) 
=& \int_1^x \phi(t)(\log t+2\gamma) dt\\
&+ \sum_{n=1}^{\infty}\tau(n)\int_1^x \phi(t)(4K_0(4\pi\sqrt{nt})-2\pi Y_0(4\pi\sqrt{nt}))dt
\end{split}\end{equation}
where $\phi$ can be any test function satisfying basic properties of regularity and $Y_0,K_0$ are Bessel functions. Equation \eqref{VoronoiSummationFormula} is the prototypical example of summation formula, for the arithmetic function $\tau(n)$. More generally, in the $GL(2)$ case, i.e. for Fourier coefficients of automoprhic forms of $GL(2)$, Voronoi summation formulas are classically well-known, see e.g. \cite[Chapter 4]{IK}.\par

Remarkably, in 2006 Miller and Schmid \cite{Miller-Schmid} proved a Voronoi formula for Maass cusp forms on $GL(3)$ twisted by additive character. Their result consists of a breakthrough and is proven by making heavy use of representation theory and the powerful method of automorphic distributions.
The main purpose of this paper is to provide a non-cuspidal version of Miller-Schmid's theorem, proving a Voronoi formula for the minimal Eisenstein series on $GL(3)$. In other words, we will derive a summation formula for the shifted triple divisor function twisted by additive character.\par

We first give some motivations for this specific choice.
Being intimately connected with the shifted correlation sums, Voronoi formulas are often encountered in the context of moments of $L$-functions \cite{DFI, SarnakRSLf}.
As a concrete example, we present the case of the moments of cusp forms, in the weight aspect.
Asymptotic formulas for moments of Hecke $L$-functions $L_f(s)$ at the central point can be conjectured by random matrix theory. More specifically, denoting by $H_k$ the normalized Hecke basis of the space of holomorphic cusp forms of weight $k$, given an integer $r$ and set of shifts $A=\{\alpha_1,\dots,\alpha_r\}$, one is interested in the shifted $r$-th moment, i.e. the quantity
$$ M_A(k) = \sum_{f\in H_k}\!\! {\vphantom{\sum}}^{h} \prod_{\alpha\in A}L_f(\tfrac{1}{2}+\alpha),$$
where the superscript $^h$ indicates the usual harmonic weight arising from Petersson norm. An asymptotic formula for $M_A(k)$ as $k\to\infty$ is known only for $r=1,2$ (see \cite{BF-Green} for the state of art results).
For bigger integers $r$, \cite{CFKRS} conjectures that
\begin{equation}\label{7june.1}
M_A(k)
\sim \sum_{V\subset A}i^{|V|k}
\left(\frac{k}{4\pi}\right)^{-2\sum_{\alpha\in V}\alpha}G_{A\smallsetminus V\cup V^-}(\tfrac{1}{2};1)
\end{equation}
where $V^-:=\{-v:v\in V\}$ and  $G_A(s;l)$ is defined by an explicit Dirichlet series, see e.g. \cite[Equation (1.3)]{CF}.
Equation \eqref{7june.1} shows the structure of the shifted $r$-moment of $L_f(\tfrac{1}{2})$; each term of the outer sum over the subsets of $A$ gives a main term, which is called an \lq\lq$\ell$-swap\rq\rq  \;term when $V$ has cardinality $\ell$. Higher swaps are classically harder to detect. For instance, the 0-swap easily comes from the diagonal term in Petersson formula, when one averages the modular coefficients. The Kloosterman sum in the off-diagonal term in Petersson formula brings into play additive characters and the so-called Estermann function
$$ D_A\bigg(s,\frac{H}{K}\bigg) = \sum_{n=1}^{\infty} \frac{\tau_A(n)e(\frac{nH}{K})}{n^s}, \quad \text{with }(H,K)=1$$
becomes crucial. Note that $\tau_A(n)$ denotes the $n$-th Dirichlet coefficient of the product $\prod_{\alpha\in A} \zeta(s+\alpha)$. The analytical properties of the Estermann function allow to detect higher swaps. For example, the poles of $D_A(s,\frac{H}{K})$ \lq\lq correspond\rq\rq \;to the 1-swap terms in \eqref{7june.1}; this is shown, under certain conditions, in \cite{CF}. One of the main tools to detect the 2-swaps is expected to be the functional equation for the Estermann function. For the cubic case, corresponding to $r=3$ and $A=\{\alpha,\beta,\gamma\}$, Proposition \ref{PropFuncEq} consists of the functional equation for the relevant Estermann function. However, obtaining the 2-swap terms by using the functional equation, and therefore deriving for example an asymptotic formula for the cubic moment, is nontrivial and requires more work, which is still in progress. \par

Let's now state our main theorem.
Given $H,K$ integers such that $(H,K)=1$, we define the \emph{shifted $GL(3)$ Estermann function}
\begin{equation}\label{EFdef}  
D_{\alpha,\beta,\gamma}\left(s,\frac{H}{K}\right) = \sum_{n=1}^{\infty}\frac{\tau_{\alpha,\beta,\gamma}(n)e(\frac{nH}{K})}{n^s},
\end{equation}
where $e(x)=e^{2\pi i x}$ and $\tau_{\alpha,\beta,\gamma}$ denotes the shifted triple divisor function, defined by 
$$\tau_{\alpha,\beta,\gamma}(n) = \sum_{abc=n} a^{-\alpha}b^{-\beta}c^{-\gamma} $$
for $\alpha,\beta,\gamma\in\mathbb C$. 

\begin{thm}\label{OurVoronoi}
For any Schwartz function $\phi$ supported in $(0,\infty)$, given $H,K$ two integers with $(H,K)=1$, denoting by $\overline H$ the inverse of $H$ modulo $K$, we have that
$$ \sum_{n=1}^\infty \tau_{\alpha,\beta,\gamma}(n)e\left(\frac{nH}{K}\right)\phi(n) $$
equals
\begin{equation}\begin{split}\notag
&\left(\frac{2\pi}{K}\right)^{\alpha+\beta+\gamma} \sum_{\varepsilon_1,\varepsilon_2,\varepsilon_3\in\{\pm 1\}}    e^{\frac{\pi}{2}i(\varepsilon_1\alpha+\varepsilon_2\beta+\varepsilon_3\gamma)} \sum_{d|K}\frac{d}{K^2}
 \sum_{h|d}\frac{\mu(h)\tau_{-\beta-\gamma,-\alpha-\gamma,-\alpha-\beta}(\frac{d}{h})}{h^{-\alpha-\beta-\gamma}}\\
&\quad  \hspace{2cm} \times
\sum_{m=1}^\infty \tau_{-\alpha,-\beta,-\gamma}(m) 
S\left(1,-\varepsilon_1\varepsilon_2\varepsilon_3\overline{H}hm;\frac{K}{d}\right)  
F_{\alpha,\beta,\gamma}\left(\frac{(2\pi)^3d^2hm}{K^3e^{-\frac{\pi}{2}i(\varepsilon_1+\varepsilon_2+\varepsilon_3)}};\phi\right)\\
&+ \tilde\phi(1-\alpha)K^{-1+\alpha}\zeta(1-\alpha+\beta)\zeta(1-\alpha+\gamma)G_{\alpha,\beta,\gamma}(1-\alpha,K) \\
&+ \tilde\phi(1-\beta)K^{-1+\beta}\zeta(1+\alpha-\beta)\zeta(1-\beta+\gamma)G_{\alpha,\beta,\gamma}(1-\beta,K) \\
&+ \tilde\phi(1-\gamma)K^{-1+\gamma}\zeta(1+\alpha-\gamma)\zeta(1+\beta-\gamma)G_{\alpha,\beta,\gamma}(1-\gamma,K),
\end{split}\end{equation}
where $S(a,b;c) = \sum_{1\leq x \leq c, (x,c)=1} e(\frac{ax+b\overline{x}}{c})$ is the Kloosterman sum, $G_{\alpha,\beta,\gamma}(s,K)$ is defined in \eqref{G3def} and
\begin{equation}\begin{split}\notag
F_{\alpha,\beta,\gamma}(x;\phi) 
&=  \frac{1}{2\pi i}\int_{(2)} \tilde\phi(1-s)   \Gamma(s-\alpha)\Gamma(s-\beta)\Gamma(s-\gamma)  x^{-s}ds\\
&= \int_{0}^{\infty} \phi(t) G_{0,3}^{3,0}(-\alpha,-\beta,-\gamma;tx) dt
\end{split}\end{equation}
with $G_{0,3}^{3,0}(-\alpha,-\beta,-\gamma;\cdot) = G_{0,3}^{3,0}([\;],[-\alpha,-\beta,-\gamma];\cdot) $ and $G_{p,q}^{m,n}([a_1,\dots,a_p],[b_1,\dots,b_q];\cdot)$ denotes the Meijer $G$-function.
\end{thm}

A first unshifted version of Theorem \ref{OurVoronoi} has been given by Ivi\'c \cite{Ivic}, then  Li \cite{Li} proved a more explicit formula. However, the proof in \cite{Li} depends on the modularity of the Eisenstein series, while our proof follows from analytical properties of the Estermann function, in particular the functional equation. Therefore, while Li uses $GL(3)$ techniques, we only rely on classical analytic number theory and some properties of the Kloosterman sums. Specializing to $\alpha=\beta=\gamma=0$ in Theorem \ref{OurVoronoi}, one recovers Li's result, except that the three residual terms collapse to one single residue of a triple pole.\par

We now compare the main formula in Theorem \ref{OurVoronoi} with that of Miller and Schmid. Theorem 1.18 in \cite{Miller-Schmid} states that, if $a_{n,m}$ are the Fourier coefficients of a cuspidal $GL(3)$ automorphic form, then for any nice test function $\phi$ we have
\begin{equation}\label{MS-AF}
\sum_{n\neq 0}a_{k,n}e(-na/c)\phi(n) 
= \sum_{d|ck}\left|\frac{c}{d}\right|\sum_{m\neq0}\frac{a_{m,d}}{|m|}S(k\overline{a},m;kc/d) \Phi\left(\frac{md^2}{c^3k}\right)
\end{equation}
where $\Phi$ denotes an explicit integral transform of $\phi$.
Theorem \ref{OurVoronoi} corresponds to the (non-cuspidal) choice $a_{1,n}$ where $a_{m,n}=a_{m,n}(\alpha,\beta,\gamma)$ is such that
$$ \sum_{m,n} \frac{a_{m,n}}{m^sn^w} = \frac{\zeta(s+\alpha+\beta)\zeta(s+\alpha+\gamma)\zeta(w+\alpha)\zeta(w+\beta)\zeta(w+\gamma)}{\zeta(s+w+\alpha+\beta+\gamma)} $$
i.e.
$$ a_{m,n} = \sum_{l|(m,n)} \mu(l)\tau_{\alpha+\beta,\alpha+\gamma,\beta+\gamma}\bigg(
\frac{m}{l}\bigg)\tau_{\alpha,\beta,\gamma}\bigg(\frac{n}{l}\bigg). $$
If we assume that $\alpha+\beta+\gamma=0$, then we can see that the right hand side of \eqref{MS-AF} corresponds to the first term of the right hand side of Theorem \ref{OurVoronoi}, being
\begin{equation}\notag
a_{m,1}=\tau_{\alpha+\beta,\alpha+\gamma,\beta+\gamma}(m)=\tau_{-\alpha,-\beta,-\gamma}(m),\quad \quad
a_{1,n}=\tau_{\alpha,\beta,\gamma}(n).
\end{equation}
Clearly, the polar terms of Theorem \ref{OurVoronoi} do not appear in the cuspidal case.\par

As usual, Theorem \ref{OurVoronoi} is stated with $\phi$ a Schwartz function to guarantee convergences, but the proof can be adapted the the case when the test function is more generic (see e.g. \cite[Lemma 4.1]{BF-Green}). For example, one can assume that $\phi:[0,\infty)\to\mathbb C$ is such that
$\tilde\phi(s)$ is regular in the region $\sigma_0<\Re(s)<\sigma_1$ and
for some $\epsilon>0$ and for all $\sigma_0<\sigma<\sigma_1$ the function $ ((1+|t|)^{-3\sigma+\frac{3}{2}+\epsilon}+1)\tilde\phi(\sigma+it) $ is integrable in $(-\infty,+\infty)$.

We remark that, speaking formally, for specific choices of the test function $\phi$ the kernel $F_{\alpha,\beta,\gamma}(x,\phi)$ can be computed explicitely. For instance, recalling the Mellin pair involving the Bessel $J$-function (see \cite{Erdelyi}, Equation (17) p. 338)
$$J_{2\nu}(2\sqrt{t})G_{0,3}^{3,0}([\;],[-\alpha,-\beta,-\gamma];xt)
\quad \longrightarrow \quad
G_{2,3}^{3,1}([1-s-\nu,1-s+\nu],[-\alpha,-\beta,-\gamma];x), $$
we have
\begin{equation}\begin{split}\notag
F_{\alpha,\beta,\gamma}(x,J_{2\nu}(2\sqrt{\cdot})) 
&= \int_{0}^{\infty} J_{2\nu}(2\sqrt{t})G_{0,3}^{3,0}([\;],[-\alpha,-\beta,-\gamma];xt) dt\\
&=G_{2,3}^{3,1}([-\nu,\nu],[-\alpha,-\beta,-\gamma];x)
\end{split}\end{equation}
and the remaining Mejier $G$-function can be expressed in terms of hypergeometric functions $\phantom{}_{2}F_2$. Similar considerations can be done for the choices $\phi(t) = K_{2\nu}(2\sqrt{t})$ the $K$-Bessel function and $\phi(t)=e^{-t}$, given the Mellin pairs (see \cite{Erdelyi}, Equation (16) and (18) p. 338)
\begin{equation}\begin{split}\notag
K_{2\nu}(2\sqrt{t})G_{0,3}^{3,0}([\;],[-\alpha,-\beta,-\gamma];xt)
&\quad \longrightarrow \quad
G_{2,3}^{3,2}([-\nu,\nu],[-\alpha,-\beta,-\gamma];x) \\
e^{-t}G_{0,3}^{3,0}([\;],[-\alpha,-\beta,-\gamma];xt)
&\quad \longrightarrow \quad
G_{1,3}^{3,1}([0],[-\alpha,-\beta,-\gamma];x).
\end{split}\end{equation}

The main theorem can be specialized to the case of the alternate sum, i.e. $H=1, K=2$, which might be of independent interest. In the case $\alpha,\beta,\gamma=0$, denoting $\tau_3$ the $3$-fold divisor function, it reads:
\begin{cor}\label{corollary_alternate}
For any Schwartz function $\phi$ supported in $(0,\infty)$, we have
\begin{equation}\begin{split}\notag
\sum_{m=1}^\infty \tau_3(m)(-1)^m\phi(m)
= \res_{s=1}&\bigg( \tilde\phi(s)\zeta(s)^3\Big( -1+\frac{6}{2^s}-\frac{6}{4^s}+\frac{2}{8^s} \Big) \bigg)\\
&\hspace{-.2cm}+ \sum_{m=1}^{\infty} \tau_3(m)\bigg( -H(m)+3H\bigg(\frac{m}{2}\bigg)-\frac{3}{2}H\bigg(\frac{m}{4}\bigg) +\frac{1}{4}H\bigg(\frac{m}{8}\bigg) \bigg),
\end{split}\end{equation}
where
\begin{equation}\notag
H(x)=\frac{1}{2\pi i}\int_{(2)}\tilde\phi(1-s)\chi(1-s)^3x^{-s}ds, \quad
\chi(1-s) = 2(2\pi)^{-s}\cos(\tfrac{\pi s}{2})\Gamma(s) .
\end{equation}
\end{cor}

We will derive this corollary from Theorem \ref{OurVoronoi}. However, it can also be proven by direct computation, noticing that
\begin{equation}\begin{split}\notag
\sum_{m=1}^\infty \tau_3(m)(-1)^m\phi(m)
= \frac{1}{2\pi i}\int_{(2)}\tilde\phi(s) \bigg(-\zeta(s)^3+2\sum_{m=1}^\infty \frac{\tau_3(2m)}{(2m)^s}\bigg)ds
\end{split}\end{equation}
and
\begin{equation}\notag
\sum_{m=1}^\infty \frac{\tau_3(2m)}{m^s} 
= \zeta(s)^3\bigg(3-\frac{3}{2^s}+\frac{1}{4^s}\bigg).
\end{equation}

\textbf{Acknowledgments}. I would like to thank Brian Conrey for suggesting the problem and for many helpful discussions. This work is supported by the FRG grant DMS 1854398 and the author is member of the INdAM group GNAMPA.

\section{The polar structure of the Estermann function}\label{SectionEstermann}

Given $H,K$ integers with $(H,K)=1$ and $A=\{\alpha_1,\dots,\alpha_r\}$ a set of shifts, we define the Estermann function
$$ D_A\left(s,\frac{H}{K}\right)=\sum_{n=1}^{\infty}\frac{\tau_A(n)e(\frac{nH}{K})}{n^s}. $$
In the following result, we analyze the polar structure of $D_A(s,\frac{H}{K})$, providing a shifted version of Conrey and Gonek's computation \cite[p. 590]{CG-HighMoments}.

\begin{lemma}\label{EF_lemma1}
The Estermann function $D_A(s,\frac{H}{K})$ admits a meromorphic continuation to the whole complex plane, with simple poles at $s=1-\alpha$, $\alpha\in A$. More precisely, $D_A(s,\frac{H}{K})$ has the same polar structure as
$$ \sum_{n=1}^\infty \frac{\tau_A(n)}{n^s}\frac{R_K(n)}{\varphi(K)} $$
where $R_K(n) := \sum_{a=1}^K{\vphantom{\sum}}' e\left( \frac{an}{K} \right)$
denotes the Ramanujan sum.
\end{lemma}

\proof
The ideas of the proof are already contained in \cite{CG-HighMoments} in the unshifted case. For the sake of completeness, here we prove the result with shifts.
For $(H,K)=1$, being 
$$ e\left( \frac{Hn}{K} \right)
=e\left(\frac{Hn/(n,K)}{K/(n,K)}\right)
=\frac{1}{\varphi(K/(n,K))}\sum_{\chi\;(K/(n,K))}\tau(\overline{\chi})\chi\left(\frac{Hn}{(n,K)}\right) ,$$
we have
\begin{equation}\notag
D_A\left(s,\frac{H}{K}\right)  
= \sum_{n=1}^{\infty} \frac{\tau_A(n)}{n^s} \frac{1}{\varphi(K/(n,K))}\sum_{\chi\;(K/(n,K))}\tau(\overline{\chi})\chi\left(\frac{Hn}{(n,K)}\right).
\end{equation}
Following Conrey and Gonek computation \cite[p. 590]{CG-HighMoments} and adding shifts $\alpha_i$, one sees that the polar part comes from the term $\chi=\chi_0$.
In particular, 
\begin{equation}\begin{split}\label{13june.0}
D_A\left(s,\frac{H}{K}\right)  
&= \sum_{d|K}\frac{1}{\varphi(\frac{K}{d})}
\sum_{\chi\;(\frac{K}{d})}\tau(\overline{\chi}) \chi(H)
\sum_{n\equiv 0\;(d)} \frac{\tau_A(n)}{n^s} \chi\left(\frac{n}{d}\right)\\
&= \sum_{d|K}\frac{d^{-s}}{\varphi(\frac{K}{d})}
\sum_{\chi\;(\frac{K}{d})}\tau(\overline{\chi}) \chi(H)
\sum_{m=1}^{\infty} \frac{\tau_A(md)\chi(m)}{m^s}.
\end{split}\end{equation}
Moreover, denoting $\nu_p(d)$ the $p$-adic evaluation of $d$, we have
\begin{equation}\begin{split}\label{13june.1}
\sum_{m=1}^{\infty} \frac{\tau_A(md)\chi(m)}{m^s}
&= \prod_{p|d} \sum_{j=0}^\infty \frac{\tau_A(p^{j+\nu_p(d)})\chi(p)^j}{p^{js}}
\times\prod_{p\nmid d}\sum_{j=0}^\infty \frac{\tau_A(p^{j})\chi(p)^j}{p^{js}}\\
&= \prod_{p|d} \frac{\sum_{j=0}^\infty \tau_A(p^{j+\nu_p(d)})\chi(p)^jp^{-js}}{\sum_{j=0}^\infty \tau_A(p^{j})\chi(p)^jp^{-js}}\times \sum_{n=1}^{\infty}\frac{\tau_A(n)\chi(n)}{n^s}.
\end{split}\end{equation}
Being 
\begin{equation}\begin{split}\label{13june.2}
\sum_{n=1}^{\infty}\frac{\tau_A(n)\chi(n)}{n^s}
= \sum_{n_1,\dots,n_r=1}^{\infty}\frac{\chi(n_1)\cdots \chi(n_r)}{n_1^{s+\alpha_1}\cdots n_r^{s+\alpha_r}}  
=\prod_{\alpha\in A}L(s+\alpha,\chi),  
\end{split}\end{equation}
this provides a meromorphic continuation of $D_A(s, \frac{H}{K})$ to the whole complex plane and shows that its only possible poles in $\Re(s) > 0$ occur at $s=1-\alpha$ and are due to the
principal character $\chi_0$ modulo $d$ for each $d$ dividing $K$. Thus, since $\tau(\overline{\chi_0}) = \mu(K/(n,K))$, the principal part of
$D_A(s, \frac{H}{K})$ is the same as that of
\begin{equation}\notag
\sum_{n=1}^{\infty} \frac{\tau_A(n)}{n^s} \frac{\mu(K/(n,K))}{\varphi(K/(n,K))}
=\sum_{n=1}^{\infty} \frac{\tau_A(n)}{n^s}  \frac{R_K(n)}{\varphi(K)} .
\end{equation}
\endproof

As mentioned in \cite[p. 591]{CG-HighMoments} a similar but more practical way to express the polar structure of $D_A(s,\frac{H}{K})$ in terms of zeta functions is given by the following result.
 
\begin{lemma}\label{EF_lemma2}
$D_A(s,\frac{H}{K})$ has the same principal part as
\begin{equation}\label{EF_PP} 
K^{-s}\prod_{\alpha\in A}\zeta(s+\alpha)G_A(s,K) 
\end{equation}
where 
\begin{equation}\label{Gdef} 
G_A(s,K) = \sum_{d|K}\frac{\mu(d)}{\varphi(d)}d^s\sum_{e|d}\frac{\mu(e)}{e^s} g_A(s,Ke/d) 
\end{equation}
and, if $K=\prod_p p^{K_p}$,
\begin{equation}\label{gdef} 
g_A(s,K) = \prod_{p|K}\left(\prod_{\alpha\in A}\left(1-\frac{1}{p^{s+\alpha}}\right)\times\sum_{j=0}^{\infty} \frac{\tau_A(p^{j+K_p})}{p^{js}}\right). 
\end{equation}
\end{lemma}

\proof
Again the bulk of the proof is already contained in \cite{CG-HighMoments}, p 591, here we address the shifted case. Let's denote $\nu_p(d)$ the $p$-adic evaluation of $d$ and 
$$ g_A(s,d,\chi) := \prod_{p|d} \frac{\sum_{j=0}^\infty \tau_A(p^{j+\nu_p(d)})\chi(p)^jp^{-js}}{\sum_{j=0}^\infty \tau_A(p^{j})\chi(p)^jp^{-js}}$$
By Equation \eqref{13june.0}, \eqref{13june.1} and \eqref{13june.2}, for $\Re(s)>1$ we have
\begin{equation}\begin{split}\notag
D_A\left(s,\frac{H}{K}\right)  
&= K^{-s}\sum_{d|K}\frac{d^s}{\varphi(d)}
\sum_{\chi\;(d)}\tau(\overline{\chi}) \chi(H)
g_A(s,K/d,\chi) \prod_{\alpha\in A}L(s+\alpha,\chi),
\end{split}\end{equation}
that provides the meromorphic continuation of $D_A(s,\frac{H}{K})$ to the whole complex plane. Similarly to what we did in the previous lemma, we look at the polar part by taking $\chi=\chi_0$
and we see that the principal part of $D_A(s,\frac{H}{K})$ is the same of that of 
\begin{equation}\begin{split}\notag
&K^{-s}\sum_{d|K}\frac{d^s}{\varphi(d)}
\tau(\overline{\chi_0}) \chi_0(H)
\sum_{m=1}^{\infty}\frac{\tau_A(m\frac{K}{d})\chi_0(m)}{m^s}.
\end{split}\end{equation}
Since $\tau(\overline{\chi_0})=\mu(d)$ and $\chi_0(m) =\sum_{e|m,e|d}\mu(e) $, the above equals
\begin{equation}\begin{split}\label{14june.1}
K^{-s}\sum_{d|K}\frac{\mu(d)}{\varphi(\frac{K}{d})}d^s \sum_{e|d}\frac{\mu(e)}{e^s}
\sum_{m=1}^{\infty}\frac{\tau_A(me\frac{K}{d})}{m^s}.
\end{split}\end{equation}
In the case of the principal character modulo 1, noticing that $g_A(s,d,\chi_0)$ corresponds to $g_A(s,d)$ defined in \eqref{gdef}, Equation \eqref{13june.1} reads
\begin{equation}\begin{split}\notag
\sum_{m=1}^{\infty} \frac{\tau_A(md)}{m^s}
= g_A(s,d) \times \prod_{\alpha\in A}\zeta(s+\alpha).
\end{split}\end{equation}
Plugging this into \eqref{14june.1}, we get
$$K^{-s}\prod_{\alpha\in A}\zeta(s+\alpha)\times\sum_{d|K}\frac{\mu(d)}{\varphi(\frac{K}{d})}d^s\sum_{e|d}\frac{\mu(e)}{e^s} g_A(s,Ke/d)$$
and \eqref{EF_PP} follows by defining $G_A(s,K)$ as in \eqref{Gdef}.
\endproof

We remark that the polar structure of $D_A(s,\tfrac{H}{K})$ does not depend on $H$.

\section{The functional equation for the $GL(3)$ Estermann function}\label{FE-EF}

In this section we specialize to the case $A=\{\alpha,\beta,\gamma\}$, and we provide a functional equation for
$$ D_{\alpha,\beta,\gamma}\bigg(s,\frac{H}{K}\bigg) 
= \sum_{n=1}^{\infty} \frac{\tau_{\alpha,\beta,\gamma}(n)e(\frac{nH}{K})}{n^s},
\hspace{.5cm} \text{with } (H,K)=1. $$
In Section \ref{SectionEstermann}, we proved that $D_{\alpha,\beta,\gamma}$ admits a meromorphic continuation to the whole complex plane, with single poles at $s=1-\alpha,1-\beta,1-\gamma$, as its principal part equals that of 
\begin{equation}\label{EF3_PP} 
K^{-s}\zeta(s+\alpha)\zeta(s+\beta)\zeta(s+\gamma)G_{\alpha,\beta,\gamma}(s,K) 
\end{equation}
where 
\begin{equation}\label{G3def} 
G_{\alpha,\beta,\gamma}(s,K) = \sum_{d|K}\frac{\mu(d)}{\varphi(d)}d^s\sum_{e|d}\frac{\mu(e)}{e^s} g_{\alpha,\beta,\gamma}(s,Ke/d) 
\end{equation}
and, if $K=\prod_p p^{K_p}$,
\begin{equation}\label{g3def}\notag 
g_{\alpha,\beta,\gamma}(s,K) = \prod_{p|K}\left(\left(1-\frac{1}{p^{s+\alpha}}\right)\left(1-\frac{1}{p^{s+\beta}}\right)\left(1-\frac{1}{p^{s+\gamma}}\right)\sum_{j=0}^{\infty} \frac{\tau_{\alpha,\beta,\gamma}(p^{j+K_p})}{p^{js}}\right). 
\end{equation}
Note that, when $p$ is a prime, the function $G_{\alpha,\beta,\gamma}(s,p)$ is conveniently easy to write; more specifically, one can prove that\footnote{Essentially this is \cite[Equation (36)]{CG-HighMoments}, with straightforward modifications (putting shifts). See also \cite{CK3}, page 7.}
$$ G_{\alpha,\beta,\gamma}(s,p) =p^s\bigg( 1 - \left(1-\frac{1}{p}\right)^{-1}\left(1-\frac{1}{p^{s+\alpha}}\right)\left(1-\frac{1}{p^{s+\beta}}\right)\left(1-\frac{1}{p^{s+\gamma}}\right) \bigg) .$$ 
In particular, at $s=1-\alpha$, we have
$$ G_{\alpha,\beta,\gamma}(1-\alpha,p) = p^{1-\alpha}\left( \frac{1}{p^{1-\alpha+\beta}} +\frac{1}{p^{1-\alpha+\gamma}} -\frac{1}{p^{2-2\alpha+\beta +\gamma}} \right) $$ 
and analogous formulas hold for $G_{\alpha,\beta,\gamma}(1-\beta,p)$ and $G_{\alpha,\beta,\gamma}(1-\delta,p)$.

Now we obtain a functional equation for $D_{\alpha,\beta,\gamma}(s,\frac{H}{K})$. 
\begin{prop}\label{PropFuncEq}
With $D_{\alpha,\beta,\gamma}(s,\frac{H}{K})$ defined as in \eqref{EFdef}, $\overline{H}H\equiv1$ (mod $K$), we have\\
Then, we have
\begin{equation}\begin{split}\notag
D_{\alpha,\beta,\gamma}\left(s,\frac{H}{K}\right)
& = K^{1-3s-\alpha-\beta-\gamma}\mathcal G_{\alpha,\beta,\gamma}(s)\sum_{\varepsilon_1,\varepsilon_2,\varepsilon_3\in\{\pm 1\}}\varepsilon_1\varepsilon_2\varepsilon_3e^{\frac{\pi}{2}i(\varepsilon_1(s+\alpha)+\varepsilon_2(s+\beta)+\varepsilon_3(s+\gamma))}\\
&\hspace{4cm} \times \sum_{d|K}\frac{1}{d^{1-2s}}\sum_{h|d}\frac{\mu(h)\tau_{-\beta-\gamma,-\alpha-\gamma,-\alpha-\beta}\left(\frac{d}{h}\right)}{h^{1-s-\alpha-\beta-\gamma}}\\ 
&\hspace{4cm} \times\sum_{a=1}^{K/d}{\vphantom{\sum}}' e\left(\frac{\overline{a}}{K/d}\right)
D_{-\alpha,-\beta,-\gamma}\left(1-s,-\varepsilon_1\varepsilon_2\varepsilon_3\frac{\overline{H}ha}{K/d} \right)
\end{split}\end{equation}
where
$$\mathcal G_{\alpha,\beta,\gamma}(s)=\mathcal G(s+\alpha)\mathcal G(s+\beta)\mathcal G(s+\gamma)$$
and
\begin{equation}\label{14june.2} \mathcal G(s) = -i(2\pi )^{s-1}\Gamma(1-s). \end{equation}
\end{prop}

Before we prove Proposition \ref{PropFuncEq}, we state a few lemmas. First, we recall the functional equation for the Hurwitz zeta-function (see e.g. \cite[Equations (14)]{Estermann} or inside the proof of Lemma 4 in \cite{Conrey2/5}).

\begin{lemma}\label{FuncEqForHurwitz}
The Hurwitz-zeta function, defined by
\begin{equation}\notag
\zeta(s,H,K) = \sum_{n\equiv H\;(K)}\frac{1}{n^s}, \quad \emph{for } \Re(s)>1 \quad(H,K\in\mathbb N)
\end{equation}
is holomorphic over $\mathbb C$ except for a simple pole at $s=1$. 
Moreover, it satisfies the functional equation
$$ \zeta(s,H,K) = K^{-s} \mathcal G(s)\left(e^{\frac{\pi}{2}is}\zeta(1-s,e(\tfrac{H}{K})) - e^{-\frac{\pi}{2}is}\zeta(1-s,e(-\tfrac{H}{K}))\right), $$
where
$$ \zeta(s,e(\tfrac{H}{K})) = \sum_{n=1}^{\infty}\frac{e(nH/K)}{n^s} $$
is the Lerch zeta-function and $\mathcal G(s)$ as in \eqref{14june.2}.
\end{lemma}

Also, denoting
\begin{equation}\label{KloostermanDef} 
S(a,b;c) = \sum_{x=1}^c{\vphantom{\sum}}'e\left(\frac{ax+b\overline{x}}{c}\right) 
\end{equation}
the usual Kloosterman sum, we prove an elementary lemma that will be useful later.

\begin{lemma}\label{LemmaExpSum}
For any $l,m,n\in\mathbb N$ and $H,K\in\mathbb Z$ with $K\neq0$ and $(H,K)=1$, $\overline{H}H=1$ modulo $K$, we have
$$ \sum_{L,M,N=1}^K e\left(\frac{HLMN+lL+mM+nN}{K}\right) = 
K\sum_{\delta|(K,m,n)}\delta\; S\left( l,-\overline{H}\frac{m}{\delta}\frac{n}{\delta};\frac{K}{\delta} \right). $$
\end{lemma}

\proof
Evaluating the sum over $N$, we get
$$  \sum_{L,M,N=1}^K e\left(\frac{HLMN+lL+mM+nN}{K}\right) = 
K \sum_{\substack{L,M=1 \\ HLM\equiv -n\;(K)}}^K e\left(\frac{lL+mM}{K}\right).$$
Writing $n=n'd$ with $d|K$ and $(n',\frac{K}{d})=1$, the above is (recall $(H,K)=1$)
\begin{equation}\begin{split}\notag 
&=K\sum_{\delta|d}\sum_{\substack{L=1 \\ (L,d)=\delta}}^K e\left(\frac{lL}{K}\right) \sum_{\substack{M=1 \\ LM\equiv -\overline{H}n'd\;(K)}}^K e\left(\frac{mM}{K}\right)\\
&= K\sum_{\delta|d}\sum_{\substack{L=1 \\ (L,\frac{d}{\delta})=1}}^{K/\delta} e\left(\frac{l\delta L}{K}\right) \sum_{\substack{M=1 \\ \delta LM\equiv -\overline{H}n'd\;(K)}}^K e\left(\frac{mM}{K}\right).
\end{split}\end{equation}
Since $\delta|d$ and $\delta|K$, the condition $\delta LM\equiv -\overline{H}n'd\;(K)$ implies $ LM\equiv -\overline{H}n'\frac{d}{\delta}\;(\frac{K}{\delta})$; moreover being $(L,\frac{d}{\delta})=1$, this yields $\frac{d}{\delta}|M$. With the change of variable $M\to\frac{d}{\delta}M$, we then obtain
$$ K\sum_{\delta|d}\sum_{\substack{L=1 \\ (L,\frac{d}{\delta})=1}}^{K/\delta} e\left(\frac{l L}{K/\delta}\right) \sum_{\substack{M=1 \\  LM\equiv -\overline{H}n'\;(\frac{K}{d})}}^{K\delta/d} e\left(\frac{mMd}{K\delta}\right) .$$
As $(n',\frac{K}{d})=1$, $LM\equiv -\overline{H}n'\;(\frac{K}{d})$ implies that $L,M$ are both coprime with $\frac{K}{d}$ and the equivalence can be written as $M\equiv -\overline{L}\overline{H}n'\;(\frac{K}{d})$. In addition, since $L$ is also coprime with with $\frac{d}{\delta}$, we also get $(L,\frac{K}{\delta})=1$, and we will denote this condition with a prime over the sum, as usual; this leads to
$$ K\sum_{\delta|d}\sum_{L=1 }^{K/\delta}{\vphantom{\sum}}' e\left(\frac{l L}{K/\delta}\right) \!\!\sum_{\substack{M=1 \\  M\equiv -\overline{L}\overline{H}n'\;(\frac{K}{d})}}^{K\delta/d} \!\!\!\!\!e\left(\frac{mMd}{K\delta}\right) 
=  K\sum_{\delta|d}\sum_{L=1 }^{K/\delta}{\vphantom{\sum}}' e\left(\frac{l L}{K/\delta}- \frac{md\overline{L}\overline{H}n'}{K\delta}\right) \sum_{r=1}^{\delta} e\left(\frac{mr}{\delta}\right)$$
since from $M\equiv -\overline{L}\overline{H}n'\;(\frac{K}{d})$, we can write $M=\frac{K}{d}r-\overline{L}\overline{H}n'$ for some $1\leq r\leq \delta$. Executing the remaining sum over $r$, the above gives
$$  K\sum_{\substack{\delta|d \\ \delta|m }}\delta\sum_{L=1 }^{K/\delta}{\vphantom{\sum}}' e\left(\frac{l L}{K/\delta}-\frac{md\overline{L}\overline{H}n'}{K\delta}\right)
=  K\sum_{\substack{\delta|d \\ \delta|m }}\delta\sum_{L=1 }^{K/\delta}{\vphantom{\sum}}' e\left(\frac{l L}{K/\delta} - \frac{\frac{m}{\delta}\frac{d}{\delta}\overline{L}\overline{H}n'}{K/\delta}\right) $$
and recalling that $n'd=n$ (and so $d|n$, also $\delta|d$ and $d|K$) we finally have
$$ K\sum_{\substack{\delta|K \\ \delta|m \\ \delta|n }}\delta\sum_{L=1 }^{K/\delta}{\vphantom{\sum}}' e\left(\frac{l L}{K/\delta}+\frac{-\frac{m}{\delta}\frac{n}{\delta}\overline{H}\overline{L}}{K/\delta}\right) 
= K\sum_{\substack{\delta|K \\ \delta|m \\ \delta|n }}\delta\;S\left( l,-\overline{H}\frac{m}{\delta}\frac{n}{\delta};\frac{K}{\delta} \right).  $$
\endproof

\begin{re}\label{8termsremark}
The above lemma can be straightforwardly straighten to the slightly more general formula
$$ \sum_{L,M,N=1}^K e\left(\frac{HLMN+\varepsilon_1 lL +\varepsilon_2 mM +\varepsilon_3 nN}{K}\right) = 
K\sum_{\delta|(K,m,n)}\delta\; S\left( l,-\varepsilon_1\varepsilon_2\varepsilon_3\overline{H}\frac{m}{\delta}\frac{n}{\delta};\frac{K}{\delta} \right) $$
where $\varepsilon_i\in\{\pm1\}$ for $i=1,2,3$.
\end{re}

We also recall Selberg's identity for Kloosterman sums, originally stated by Selberg \cite{Selberg} without a proof, then proved by Kuznetsov \cite{Kuznetsov}.

\begin{lemma}\label{SelbergIdentity}
For any positive integer $q$ and for any integers $a,b$, we have
$$ S(a,b;q) = \sum_{g|(a,b,q)}g\; S\left(1,\frac{ab}{g^2};\frac{q}{g}\right). $$
\end{lemma}

We are now ready to prove Proposition \ref{PropFuncEq}.

\proof
By definition of $\tau_{\alpha,\beta,\gamma}$, splitting the sums over residue classes (mod $K$), we have
\begin{equation}\label{3mar.3}
D_{\alpha,\beta,\gamma}\left(s,\frac{H}{K}\right)
=\sum_{L,M,N=1}^{K}e\left(\frac{HLMN}{K}\right)\zeta(s+\alpha,L,K)\zeta(s+\beta,M,K)\zeta(s+\gamma,N,K).
\end{equation}
Now we apply three times the functional equation for the Hurwitz zeta-function (see Lemma \ref{FuncEqForHurwitz}).
Doing so, we write the quantity in \eqref{3mar.3} as a sum of 8 terms, the first being 
\begin{equation}\begin{split}\label{25july.2}
\frac{e^{\frac{\pi}{2}i(3s+\alpha+\beta+\gamma)}}{K^{3s+\alpha+\beta+\gamma}}\mathcal G_{\alpha,\beta,\gamma}(s)&\sum_{L,M,N=1}^K e\left(\frac{HLMN}{K}\right) \\
&\times\zeta(1-s-\alpha,e(\tfrac{L}{K}))\zeta(1-s-\beta,e(\tfrac{M}{K}))\zeta(1-s-\gamma,e(\tfrac{N}{K}))
\end{split}\end{equation}
where $\mathcal G_{\alpha,\beta,\gamma}(s)=\mathcal G(s+\alpha)\mathcal G(s+\beta)\mathcal G(s+\gamma)$. 
The above is
$$\frac{e^{\frac{\pi}{2}i(3s+\alpha+\beta+\gamma)}}{K^{3s+\alpha+\beta+\gamma}}\mathcal G_{\alpha,\beta,\gamma}(s)\sum_{l,m,n}\frac{l^\alpha m^\beta n^\gamma}{(lmn)^{1-s}}\sum_{L,M,N=1}^K e\left(\frac{HLMN+lL+mM+nN}{K}\right) .$$
We focus on the sums above, ignoring for now the factor in front.
The inner exponential sum can be evaluated by Lemma \ref{LemmaExpSum}, getting
\begin{equation}\begin{split}\notag
&K\sum_{l,m,n}\frac{l^\alpha m^\beta n^\gamma}{(lmn)^{1-s}}\sum_{\delta|(K,m,n)}\delta \;S\left(l,-\overline{H}\frac{m}{\delta}\frac{n}{\delta};\frac{K}{\delta}\right)\\
=& K\sum_{\delta|K}\frac{\delta^{\beta+\gamma}}{\delta^{1-2s}}\sum_{l,m,n}\frac{l^\alpha m^\beta n^\gamma}{(lmn)^{1-s}}S\left(l,-\overline{H}mn;\frac{K}{\delta}\right)\\
=& K\sum_{\delta|K}\frac{\delta^{\beta+\gamma}}{\delta^{1-2s}}\sum_{l,r}\frac{l^\alpha \tau_{-\beta,-\gamma}(r)}{(lr)^{1-s}}S\left(l,-\overline{H}r;\frac{K}{\delta}\right)  
\end{split}\end{equation}
and an application of Selberg's identity (see Lemma \ref{SelbergIdentity} and remember that $(H,K)=1$) yields
\begin{equation}\begin{split}\label{3mar.4} 
&K\sum_{\delta|K}\frac{\delta^{\beta+\gamma}}{\delta^{1-2s}}\sum_{l,r}\frac{l^\alpha \tau_{-\beta,-\gamma}(r)}{(lr)^{1-s}} \sum_{g|(l,r,\frac{K}{\delta})}g\;S\left(1,-\overline{H}\frac{lr}{g^2};\frac{H}{\delta g}\right)  \\
=& K\sum_{\delta|K}\sum_{g|\frac{K}{\delta}}\frac{\delta^{\beta+\gamma}}{\delta^{1-2s}} \frac{g^\alpha}{g^{1-2s}} \sum_{l,r}\frac{l^\alpha \tau_{-\beta,-\gamma}(rg)}{(lr)^{1-s}} S\left(1,-\overline{H}lr;\frac{K}{\delta g}\right) .
\end{split}\end{equation}
The well-known multiplicativity relation for the divisor function:
$$ \tau_{-\beta,-\gamma}(rg) = \sum_{\substack{h|r \\ h|g}} \frac{\mu(h)}{h^{-\beta-\gamma}} \tau_{-\beta,-\gamma}\left(\frac{r}{h}\right) \tau_{-\beta,-\gamma}\left(\frac{g}{h}\right)$$
allows us to rewrite \eqref{3mar.4} as 
\begin{equation}\begin{split}\notag
& K\sum_{\delta|K}\sum_{g|\frac{K}{\delta}}\frac{\delta^{\beta+\gamma}g^\alpha}{(\delta g)^{1-2s}} \sum_{h|g} \frac{\mu(h)}{h^{1-s-\beta-\gamma}}\tau_{-\beta,-\gamma}\left(\frac{g}{h}\right) \sum_{l,r}\frac{l^{\alpha} \tau_{-\beta,-\gamma}(r)}{(lr)^{1-s}} S\left(1,-\overline{H}hlr;\frac{K}{\delta g}\right) \\
=& K\sum_{d|K}\frac{1}{d^{1-2s}} \sum_{d=\delta g}\delta^{\beta+\gamma}g^{\alpha} \sum_{h|g} \frac{\mu(h)}{h^{1-s-\beta-\gamma}}\tau_{-\beta,-\gamma}\left(\frac{g}{h}\right) \sum_{l,m,n}\frac{l^{\alpha}m^{\beta}n^{\gamma}}{(lmn)^{1-s}} S\left(1,-\overline{H}hlmn;\frac{K}{d}\right) \\
=& K\sum_{d|K}\frac{1}{d^{1-2s}}\sum_{h|d}\frac{\mu(h)}{h^{1-s-\alpha-\beta-\gamma}}\bigg[ \sum_{d=\delta g h}\delta^{\beta+\gamma}g^{\alpha}\tau_{-\beta,-\gamma}(g)\bigg] \sum_{n=1}^{\infty}\frac{\tau_{-\alpha,-\beta,-\gamma}(n)}{n^{1-s}} S\left(1,-\overline{H}hn;\frac{K}{d}\right) .
 \end{split}\end{equation}
Finally, the remaining bracket is:
\begin{equation}\begin{split}\notag
\sum_{d=\delta g h}\delta^{\beta+\gamma}g^{\alpha}\tau_{-\beta,-\gamma}(g) 
&=\sum_{d=\delta ab h}\delta^{\beta+\gamma}(ab)^{\alpha}a^{\beta}b^{\gamma} \\
&=\sum_{\frac{d}{h}=\delta ab}\delta^{\beta+\gamma}a^{\alpha+\beta}b^{\alpha+\gamma}
=\tau_{-\beta-\gamma,-\alpha-\gamma,-\alpha-\beta}\bigg(\frac{d}{h}\bigg) .
\end{split}\end{equation}
Therefore, the term in \eqref{25july.2} gives
\begin{equation}\begin{split}\notag
&\frac{e^{\frac{\pi}{2}i(3s+\alpha+\beta+\gamma)}}{K^{-1+3s+\alpha+\beta+\gamma}}\mathcal G_{\alpha,\beta,\gamma}(s)
\sum_{d|K}\frac{1}{d^{1-2s}}
\sum_{h|d}\frac{\mu(h)}{h^{1-s-\alpha-\beta-\gamma}}\tau_{-\beta-\gamma,-\alpha-\gamma,-\alpha-\beta}\left(\frac{d}{h}\right)\\
&\hspace{7cm} \times\sum_{n=1}^{\infty}\frac{\tau_{-\alpha,-\beta,-\gamma}(n)}{n^{1-s}} S\left(1,-\overline{H}hn;\frac{K}{d}\right) \\
=& \frac{e^{\frac{\pi}{2}i(3s+\alpha+\beta+\gamma)}\mathcal G_{\alpha,\beta,\gamma}(s)}{K^{-1+3s+\alpha+\beta+\gamma}}
\sum_{d|K}\frac{1}{d^{1-2s}}
\sum_{h|d}\frac{\mu(h)\tau_{-\beta-\gamma,-\alpha-\gamma,-\alpha-\beta}\left(\frac{d}{h}\right)}{h^{1-s-\alpha-\beta-\gamma}} \\
&\hspace{7cm} \times \sum_{a=1}^{K/d}{\vphantom{\sum}}' e\left(\frac{\overline{a}}{K/d}\right) 
D_{-\alpha,-\beta,-\gamma}\left(1-s, -\frac{\overline{H}ha}{K/d} \right),
\end{split}\end{equation}
by definition of the Kloosterman sum \eqref{KloostermanDef}.
The above calculation concludes the proof, since the other terms can be treated analogously in view of Remark \ref{8termsremark}.
\endproof

\section{Proof of Theorem \ref{OurVoronoi}}\label{sectionVoronoi}

Now that we studied the polar structure and obtained the functional equation for the Estermann function $D_{\alpha,\beta,\gamma}(s,\frac{H}{K})$, we are ready to prove the Voronoi summation formula, i.e. Theorem \ref{OurVoronoi}.
The strategy is pretty classical; using Perron formula and performing a contour shift, in view of \eqref{EF3_PP}, we get
\begin{equation}\begin{split}\notag
\sum_{n=1}^\infty &\tau_{\alpha,\beta,\gamma}(n)e\left(\frac{nH}{K}\right)\phi(n)
= \frac{1}{2\pi i}\int_{(2)} \tilde\phi(s)D_{\alpha,\beta,\gamma}\left(s,\frac{H}{K}\right)ds \\
&= \frac{1}{2\pi i}\int_{(-1)} \tilde\phi(s)D_{\alpha,\beta,\gamma}\left(s,\frac{H}{K}\right)ds + \bigg(\res_{s=1-\alpha}+\res_{s=1-\beta}+\res_{s=1-\gamma}\bigg)\left( \tilde\phi(s)D_{\alpha,\beta,\gamma}\left(s,\frac{H}{K}\right) \right)\\
&= \mathcal I_{\alpha,\beta,\gamma}(H,K) 
+ \tilde\phi(1-\alpha)K^{-1+\alpha}\zeta(1-\alpha+\beta)\zeta(1-\alpha+\gamma)G_{\alpha,\beta,\gamma}(1-\alpha,K)\\
&\hspace{2.75cm}+ \tilde\phi(1-\beta)K^{-1+\beta}\zeta(1+\alpha-\beta)\zeta(1-\beta+\gamma)G_{\alpha,\beta,\gamma}(1-\beta,K) \\
&\hspace{2.75cm}+ \tilde\phi(1-\gamma)K^{-1+\gamma}\zeta(1+\alpha-\gamma)\zeta(1+\beta-\gamma)G_{\alpha,\beta,\gamma}(1-\gamma,K),
\end{split}\end{equation}
with
\begin{equation}\begin{split}\notag
\mathcal I = \mathcal I_{\alpha,\beta,\gamma}(H,K) 
:&= \frac{1}{2\pi i}\int_{(-1)} \tilde\phi(s)D_{\alpha,\beta,\gamma}\left(s,\frac{H}{K}\right)ds \\
&= \frac{1}{2\pi i}\int_{(2)} \tilde\phi(1-s)D_{\alpha,\beta,\gamma}\left(1-s,\frac{H}{K}\right)ds.
\end{split}\end{equation}
Applying the functional equation for the Estermann function given by Proposition \ref{PropFuncEq}, the above can be written as 
\begin{equation}\begin{split}\notag
\mathcal I
&= \frac{1}{2\pi i}\int_{(2)} \tilde\phi(1-s) \mathcal G(1-s+\alpha)\mathcal G(1-s+\beta)\mathcal G(1-s+\gamma)K^{-2+3s-\alpha-\beta-\gamma} \\
&\quad  \times   \sum_{\varepsilon_1,\varepsilon_2,\varepsilon_3\in\{\pm 1\}}  \varepsilon_1\varepsilon_2\varepsilon_3 e^{\frac{\pi}{2}i[\varepsilon_1(1-s+\alpha)+\varepsilon_2(1-s+\beta)+\varepsilon_3(1-s+\gamma)]}\sum_{d|K}\frac{1}{d^{2s-1}}\\
&\quad  \times \sum_{h|d}\frac{\mu(h)\tau_{-\beta-\gamma,-\alpha-\gamma,-\alpha-\beta}(\frac{d}{h})}{h^{s-\alpha-\beta-\gamma}} \sum_{a=1}^{K/d}{\vphantom{\sum}}' e\left(\frac{\overline{a}}{K/d}\right) 
\sum_{m=1}^\infty \frac{\tau_{-\alpha,-\beta,-\gamma}(m)}{m^{s}}e\left(-\varepsilon_1\varepsilon_2\varepsilon_3\frac{\overline{H}ham}{K/d}\right) ds \\
&= \sum_{\varepsilon_1,\varepsilon_2,\varepsilon_3\in\{\pm 1\}}  \varepsilon_1\varepsilon_2\varepsilon_3 i  e^{\frac{\pi}{2}i(\varepsilon_1+\varepsilon_2+\varepsilon_3)}  e^{\frac{\pi}{2}i(\varepsilon_1\alpha+\varepsilon_2\beta+\varepsilon_3\gamma)}  (2\pi)^{\alpha+\beta+\gamma} K^{-2-\alpha-\beta-\gamma} \\
&\quad  \times \sum_{d|K}d  \sum_{h|d}\frac{\mu(h)\tau_{-\beta-\gamma,-\alpha-\gamma,-\alpha-\beta}(\frac{d}{h})}{h^{-\alpha-\beta-\gamma}} \sum_{m=1}^\infty \tau_{-\alpha,-\beta,-\gamma}(m)\sum_{a=1}^{K/d}{\vphantom{\sum}}' e\left(\frac{\overline{a}}{K/d}-\varepsilon_1\varepsilon_2\varepsilon_3\frac{\overline{H}ham}{K/d}\right)\\
&\quad  \times \frac{1}{2\pi i}\int_{(2)} \tilde\phi(1-s) (2\pi)^{-3s} \Gamma(s-\alpha)\Gamma(s-\beta)\Gamma(s-\gamma)\left(\frac{K^{3}}{d^{2}hm}\right)^s e^{-\frac{\pi}{2}is(\varepsilon_1+\varepsilon_2+\varepsilon_3)} ds.
\end{split}\end{equation}
Note that $  \varepsilon_1\varepsilon_2\varepsilon_3 i  e^{\frac{\pi}{2}i(\varepsilon_1+\varepsilon_2+\varepsilon_3)} =1 $ for all $\varepsilon_1,\varepsilon_2,\varepsilon_3\in\{\pm 1\}$.
Therefore, defining
\begin{equation}\begin{split}\notag
F_{\alpha,\beta,\gamma}(x,\phi) 
:&= \frac{1}{2\pi i}\int_{(2)} \tilde\phi(1-s) \Gamma(s-\alpha)\Gamma(s-\beta)\Gamma(s-\gamma)x^{-s}  ds \\
\end{split}\end{equation}
we have
\begin{equation}\begin{split}\notag
\mathcal I
&= \frac{1}{K^2}\bigg(\frac{2\pi}{K}\bigg)^{\alpha+\beta+\gamma}\sum_{\varepsilon_1,\varepsilon_2,\varepsilon_3\in\{\pm 1\}}  e^{\frac{\pi}{2}i(\varepsilon_1\alpha+\varepsilon_2\beta+\varepsilon_3\gamma)} \sum_{d|K}d  \sum_{h|d}\frac{\mu(h)\tau_{-\beta-\gamma,-\alpha-\gamma,-\alpha-\beta}(\frac{d}{h})}{h^{-\alpha-\beta-\gamma}}\\
&\quad\quad\quad\quad\quad  \times \sum_{m=1}^\infty \tau_{-\alpha,-\beta,-\gamma}(m) S\left(1,-\varepsilon_1\varepsilon_2\varepsilon_3\overline{H}hm;\frac{K}{d}\right)  F_{\alpha,\beta,\gamma}\left(\frac{(2\pi)^{3}d^{2}hm}{K^{3}e^{-\frac{\pi}{2}i(\varepsilon_1+\varepsilon_2+\varepsilon_3)}};\phi\right).
\end{split}\end{equation}

Writing $F_{\alpha,\beta,\gamma}(x;\phi)$ in terms of the Meijer $G$-function is essentially straightforward, as by Mellin inversion we get
\begin{equation}\begin{split}\notag
F_{\alpha,\beta,\gamma}(x;\phi)
&= \frac{1}{2\pi i}\int_{(2)}\tilde\phi(1-s) \Gamma(s-\alpha)\Gamma(s-\beta)\Gamma(s-\gamma) x^{-s} ds \\
&= \int_{0}^{\infty} \phi(t) \frac{1}{2\pi i}\int_{(2)} \Gamma(s-\alpha)\Gamma(s-\beta)\Gamma(s-\gamma)(tx)^{-s}ds\;dt \\
&= \int_{0}^{\infty} \phi(t) G_{0,3}^{3,0}([\;],[-\alpha,-\beta,-\gamma];tx) dt
\end{split}\end{equation}
since we know that (see \cite[p 353 Equation (49)]{Erdelyi}), denoting by $G_{p,q}^{m,n}$ the Meijer $G$-function
$$ G_{p,q}^{m,n}([a_1,\dots,a_p],[b_1,\dots,b_q];y) = \frac{1}{2\pi i}\int_{L}\frac{\prod_{j=1}^{m}\Gamma(s+b_j)\prod_{j=1}^{n}\Gamma(1-s-a_j)}{\prod_{j=m+1}^{q}\Gamma(1-s-b_j)\prod_{j=n+1}^{p}\Gamma(s+a_j)}y^{-s}ds $$
where $L$ is a path separating the poles of $\Gamma(s+b_j)$ and those of $\Gamma(1-s-a_j)$.

\section{Proof of Corollary \ref{corollary_alternate}}

Let's specialize Theorem \ref{OurVoronoi} to the case where $\alpha=\beta=\gamma=0$ and $H=1, K=2$. 
We recall the functional equation for the Riemann zeta function, being
\begin{equation}\notag  \zeta(1-s) = \chi(1-s)\zeta(s) \end{equation}
where
\begin{equation}\begin{split}\notag
\chi(1-s) = 2(2\pi)^{-s}\cos(\tfrac{\pi s}{2})\Gamma(s) 
= (2\pi)^{-s}\Gamma(s)\Big(e^{\frac{\pi i s}{2}} + e^{-\frac{\pi i s}{2}}\Big)
= \pi^{\frac{1}{2}-s}\frac{\Gamma(\frac{s}{2})}{\Gamma(\frac{1-s}{2})},
\end{split}\end{equation}
with
$$ \chi(1-s) = \bigg(\frac{|t|}{2\pi}\bigg)^{s-1}\Big(1+O(|t|^{-1})\Big). $$

When the shifts $\alpha,\beta,\gamma$ all equal 0 and $H=1,K=2$, Theorem \ref{OurVoronoi} reads
\begin{equation}\notag
 \sum_{m=1}^\infty \tau_3(m)(-1)^m\phi(m) =\mathcal R + \mathcal I
\end{equation}
where
$$ \mathcal R = \res_{s=1}\bigg( \tilde\phi(s)\sum_{m=1}^{\infty}\frac{d_3(m)e(\frac{m}{2})}{m^s} \bigg)
= \res_{s=1}\bigg( \tilde\phi(s)\zeta(s)^3\bigg( -1+\frac{6}{2^s}-\frac{6}{4^s}+\frac{2}{8^s}  \bigg)\bigg)$$
and
\begin{equation}\begin{split}\label{26aug_2}
\mathcal I =&
\sum_{\varepsilon_1,\varepsilon_2,\varepsilon_3\in\{\pm 1\}}  
\sum_{d|2}\frac{d}{4}
\sum_{h|d}\mu(h)\tau_3\bigg(\frac{d}{h}\bigg)\\
&\hspace{2cm} \times
\sum_{m=1}^\infty \tau_3(m) 
S\left(1,-\varepsilon_1\varepsilon_2\varepsilon_3hm;\frac{2}{d}\right)  
F\left(\frac{(2\pi)^3d^2hm}{2^3e^{-\frac{\pi}{2}i(\varepsilon_1+\varepsilon_2+\varepsilon_3)}}\right) 
\end{split}\end{equation}
with
\begin{equation}\notag
F(x) = F_{0,0,0}(x;\phi) 
=  \frac{1}{2\pi i}\int_{(2)} \tilde\phi(1-s)   \Gamma(s)^3  x^{-s}ds.
\end{equation}
Expanding out the three terms ($d=1,h=1$ and $d=2,h=1$ and $d=2,h=2$) in \eqref{26aug_2},
noticing that $S(1,\pm m;1)=1$ and
\begin{equation}\notag
S(1, m;2) = e\left( \frac{m+1}{2} \right) = (-1)^{m+1} = S(1,-m;2),
\end{equation}
we have
\begin{equation}\begin{split}\notag
\mathcal I 
&=\frac{1}{4}
\sum_{\varepsilon_1,\varepsilon_2,\varepsilon_3\in\{\pm 1\}} 
\sum_{m=1}^\infty \tau_3(m) \bigg[
(-1)^{m+1}  F\left(\frac{\pi^3m}{e^{-\frac{\pi}{2}i(\varepsilon_1+\varepsilon_2+\varepsilon_3)}}\right)\\
&\hspace{3.5cm}+ 6  F\left(\frac{4\pi^3m}{e^{-\frac{\pi}{2}i(\varepsilon_1+\varepsilon_2+\varepsilon_3)}}\right)
- 2 F\left(\frac{8\pi^3m}{e^{-\frac{\pi}{2}i(\varepsilon_1+\varepsilon_2+\varepsilon_3)}}\right)
\bigg]\\
&=\frac{1}{4} \sum_{m=1}^\infty \tau_3(m) \frac{1}{2\pi i}\int_{(2)}\tilde\phi(1-s)\Gamma(s)^3 (2\pi)^{-3s}
\sum_{\varepsilon_1,\varepsilon_2,\varepsilon_3\in\{\pm 1\}} e^{s\frac{\pi}{2}i(\varepsilon_1+\varepsilon_2+\varepsilon_3)}\\
&\hspace{3.5cm} \times\bigg[
(-1)^{m+1}  \left(\frac{8}{m}\right)^{s}+ 6  \left(\frac{2}{m}\right)^{s}
- 2 \left(\frac{1}{m}\right)^{s}\bigg].
\end{split}\end{equation}
Since
$$ \sum_{\varepsilon_1,\varepsilon_2,\varepsilon_3\in\{\pm 1\}} e^{-\frac{\pi}{2}is(\varepsilon_1+\varepsilon_2+\varepsilon_3)} =  8\cos^3\left(\frac{\pi s}{2}\right), $$
denoting $\chi(1-s) = 2\cos(\frac{\pi s}{2})\Gamma(s)(2\pi)^{-s}$,
we get
\begin{equation}\begin{split}\label{18june.1}
\mathcal I 
& =\sum_{m=1}^\infty \tau_3(m)
 \bigg[
\frac{(-1)^{m+1}}{4} H\left(\frac{m}{8}\right)
+ \frac{3}{2} H\left(\frac{m}{2}\right)
- \frac{1}{2} H\left(m\right)
\bigg].
\end{split}\end{equation}
We notice that
\begin{equation}\begin{split}\notag
\sum_{m=1}^\infty \frac{\tau_3(2m)/\tau_3(2)}{m^s} 
&= \prod_{p\neq 2}\sum_{j=0}^\infty \frac{\tau_3(2p^j)/\tau_3(2)}{p^{js}} \times \sum_{j=0}^\infty \frac{\tau_3(2^{j+1})/\tau_3(2)}{2^{js}}\\
&= \zeta(s)^3 \times \bigg(1-\frac{1}{2^s}\bigg)^3 \times \frac{\frac{1}{3}(3-\frac{3}{2^s}+\frac{1}{4^s})}{(1-2^{-s})^3},
\end{split}\end{equation}
being 
\begin{equation}\begin{split}\notag
\sum_{j=0}^\infty \frac{\tau_3(2^{j+1})}{\tau_3(2)}X^j 
= \sum_{j=0}^\infty \frac{1}{\tau_3(2)}\binom{j+3}{2}X^j
=  \frac{1}{\tau_3(2)}\frac{3-3X+X^2}{(1-X)^3}.
\end{split}\end{equation}
Therefore
\begin{equation}\label{16june.1}
\sum_{m=1}^\infty \frac{\tau_3(2m)}{m^s} 
= \zeta(s)^3\bigg(3-\frac{3}{2^s}+\frac{1}{4^s}\bigg). 
\end{equation}
Hence, the first term in \eqref{18june.1} can be expanded as follows
\begin{equation}\begin{split}\label{18june.2}
\sum_{m=1}^\infty \tau_3(m) \frac{(-1)^{m+1}}{4} H\left(\frac{m}{8}\right)
&=\sum_{m=1}^\infty \tau_3(m) \frac{1}{4} H\left(\frac{m}{8}\right) - 2\sum_{m \text{ even}} \tau_3(m) \frac{1}{4} H\left(\frac{m}{8}\right)\\
&= \frac{1}{4}\sum_{m=1}^\infty \tau_3(m) H\left(\frac{m}{8}\right) - \frac{1}{2}\sum_{m=1}^\infty \tau_3(2m) H\left(\frac{m}{4}\right)\\
&= \sum_{m=1}^\infty \tau_3(m)\bigg( \frac{1}{4}H\left(\frac{m}{8}\right) - \frac{3}{2}H\left(\frac{m}{4}\right)+\frac{3}{2} H\left(\frac{m}{2}\right)-\frac{1}{2}H(m)\bigg),
\end{split}\end{equation}
since
\begin{equation}\begin{split}\notag
\sum_{m=1}^\infty \tau_3(2m) H\left(\frac{m}{4}\right) 
&=\sum_{m=1}^\infty \tau_3(2m) \frac{1}{2\pi i}\int_{(2)} \tilde\phi(1-s)\chi(1-s)^3\left(\frac{m}{4}\right)^{-s}ds \\
&= \frac{1}{2\pi i}\int_{(2)} \tilde\phi(1-s)\chi(1-s)^3 4^s \zeta(s)^3\bigg(3-\frac{3}{2^s}+\frac{1}{4^s}\bigg)ds \\
&=\sum_{m=1}^\infty \tau_3(m)   \bigg(3H\left(\frac{m}{4}\right)-3 H\left(\frac{m}{2}\right)+H(m)\bigg)
\end{split}\end{equation}
by \eqref{16june.1}. 
The claim follows plugging \eqref{18june.2} into \eqref{18june.1}.
\endproof

{\small

\end{document}